\documentclass[10pt]{amsart}
\usepackage{amsmath,amsthm,amscd,amsfonts,amssymb}
\usepackage[all]{xy}

\begin{document}
   \baselineskip   .7cm
\newtheorem{defn}{Definition}[section]
\theoremstyle{definition}
\newtheorem{thm}[defn]{Theorem}
\newtheorem{propo}[defn]{Proposition}
\newtheorem{cor}[defn]{Corollary}
\newtheorem{lem}[defn]{Lemma}
\theoremstyle{remark}
\newtheorem{rem}{Remark}
 \newtheorem{ex}[defn]{Example}
\renewcommand\o{{{\mathcal O}}}
\newcommand\s{\sigma}
\newcommand\w{\widehat}
\newcommand\cal{\mathcal}
\newcommand\limi{\varinjlim}
\newcommand\limp{\varprojlim}
\newcommand\te{\text}

\title [Ideal class group anhilators]{Ideal class group annihilators}

\author[A. \'Alvarez]{A. \'Alvarez${}^*$     }
 \dedicatory{A Jes\'us Mu\~noz D{\'{\i}}az por su   cumplea\~nos}
\address{Departamento de Matem\'aticas  \\
Universidad de Salamanca \\ Plaza de la Merced 1-4. Salamanca
(37008). Spain.}

\thanks{$*$  Departamento de Matem\'aticas.
Universidad de Salamanca. Spain.}

\maketitle

 \begin{abstract}
 We study certain correspondences over Drinfeld modular varieties  given by   sums of Hecke correspondences.
  We propose generalizations of Stickelberger's theorem for
 higher dimensions. Using this result, we study anihilators for some cusp
 forms.
 \end{abstract}

\tableofcontents


\section{Introduction.}
The aim of this article is to  propose generalizations of
Stickelberger's theorem for
 higher dimensions. Using these results, we study anihilators for some cusp
 forms. We address  certain correspondences, given
by sums of Hecke correspondences and defined over Drinfeld modular
varieties .

Let ${\Bbb P}^1$ be the projective line scheme over ${\Bbb F}_q$,
${\Bbb P}^1\setminus \{\infty\}=Spec({\Bbb F}_q[t])$ and let
  $   I=p(t){\Bbb F}_q[t]$   be an   ideal   in ${\Bbb F}_q[t]$ with
$deg(p(t))= {d+1}$.
  There exists an
abelian Galois extension, $K^\infty_I/ {\Bbb F}_q(t)$, of group
$G_{ I}\simeq ({\Bbb F}_q[t]/I)^\times$. These fields are the
Carlitz extensions and   are the cyclotomic fields in the  case of
function fields. c.f. \cite{Ca}.

Let us consider    the $S$-incomplete $L$-function evaluator ($ S:
=\vert I\vert\cup \{\infty\}$)
 $$ \underset{x\in \vert {\Bbb P}^1\vert \setminus S  }\prod (1-\tau_x\cdot z^{deg(x)})^{-1},  $$
 $\tau_x\in G_{  I}$ being the Frobenius element for $x\in \vert {\Bbb P}^1\vert$. This Euler product can be expressed
 as:
$$Q(z)+\frac{ (\sum_{h\in G_{  I}}h)\cdot z^{ d+1}} {1-q\cdot z}, $$
 $Q(z)$ being    a polynomial  in $ {\Bbb
Z}[G_{ I} ][z]$ of degree $ d$. If  one denotes $
Q(z):=\sum_{i=0}^{ d}\gamma_i\cdot z^i$,  with     $\gamma_i\in
{\Bbb Z}[G_{ I} ]$, then the correspondence
$$ \sum_{i=0}^{ d}\Gamma(Fr^{ d-i})* \Gamma(\gamma_i)$$
is trivial on $Spec(K^\infty_I \otimes K^\infty_I)$. This is
proved  for $S=\{0,1, \infty\}$   in {\cite C} and for general
$S=\{?, \infty\}$   in \cite{An1}.  This result is   analogous to
the function field case of   Stickelberger's theorem. Here,
$\Gamma(Fr^i)$ denotes the graphic of the Frobenius morphism,
$Fr^i$, and $\Gamma(\gamma_i)$ is a sum of graphics of elements of
$G_{  I}$. For arbitrary smooth curves   analogous   results can
be found in \cite{Al2}.

   These trivial correspondences give  an annihilating polynomial for the
operator given by  the correspondence $\Gamma(Fr )$ acting on the
${\Bbb Q}[G_{  I}]$-module, $H^1((Y^\infty_{ I})_{\bar {\Bbb F} },
{\Bbb Q}_l)$,  and this implies proofs of the Brumer-Stark
conjecture in the function field case  (\cite{An1}, \cite{C},
\cite{H1}, \cite{Ta}, \cite{Al2}). $Y^\infty_{  I}$ denotes the
Riemann variety associated with $K^\infty_{ I}/{\Bbb F}_q$.

Here, we study   the Euler products
 $$\prod_{x\in \vert {\Bbb P}^1\vert \setminus S  }\frac{1}{ 1-\sigma^x_1\cdot z +  q \sigma^x_2\cdot z^2
-   \cdots +              (-1)^n  q^{n(n-1)/2 } \sigma^x_n\cdot
z^n }
   =\sum_{m \geq 0} T(m)\cdot z^m,   $$

 where, "$T(m)$"  and "$\sigma^x_j$" are Hecke correspondences over certain modular Drinfeld varieties of dimension $n$,
 ${\mathcal E}^{I\infty}_{n, ?}$.
For   precision in the notation,    see section 2.2.

We prove

{\bf {Theorem 1}}
 The correspondence
 $$T(n\cdot d)+\Gamma(Fr)* T(n\cdot d-1)+\cdots +\Gamma(Fr^{n\cdot d-1})* T (1)+\Gamma(Fr^{n\cdot d }) $$
is trivial(=rationally equivalent to $0$ as cycles) in
 ${\mathcal
E}^{I\infty}_{n, ?} \times {\mathcal E}^{I\infty}_{n, ?}.$

$*$ denotes the product of correspondences. This result for $n=1$
gives us Stickelberger's theorem for the cyclotomic function
fields, \cite{An1}. To do so, we study the isogenies of Drinfeld
modules, as stated in \cite{Gr2}.

${\mathcal E}^{I\infty}_{n, ?}$ are affine schemes over
$Spec({\Bbb F}_q[t,1/h(t)]$,  and  $h(t)$ is a polynomial which
depends on $I$. We denote $   {\mathcal E}_2({I\infty}):={\mathcal
E}^{I\infty}_{2, ?}\otimes_{{\Bbb F}_q[t ]}{\Bbb F}_q(t)$.
  ${\mathcal E}({I\infty})$ is a
smooth affine curve which is defined over $K^\infty_I$. $\overline
{{\mathcal E}}({I\infty}) $ denotes the  projective curve over
$K_I$ associated with ${\mathcal E}({I\infty})$. Theorem 1 has the
following results as a sequel.

{\bf {Lemma}}  $T(2d )+  T(2d -1)+\cdots +  T (1)+\Gamma(Id) $
annihilates the group $Pic({\mathcal E}({I\infty}))$. It seems to
be a Stickelberger's  theorem for the affine modular curve
${\mathcal E}({I\infty})$   over $K^\infty_I$.

There exists an arithmetic subgroup, $\Gamma_{I\infty}$, of
$Gl_2({\Bbb F}_q[t])$  such that if we denote by $\Omega$ the
Drinfeld's upper half-plane and by $\overline {M}_{
\Gamma_{I\infty} }$   the smooth projective model of the algebraic
curve associated with $\Omega/\Gamma_{I\infty}$, we have:
$$\overline {M}_{ \Gamma_{I\infty} }=\overline {\mathcal
E}({I\infty})\otimes_{K^\infty_I}C, $$    $C$  being   the
algebraic closure of the completion of ${\Bbb F}_q(t)$ at
$\infty$. In the usual way cusp forms for  $\Gamma_{I\infty}$ are
given by $H^0(\overline {M}_{ \Gamma_{I\infty} },
\Omega_{\overline {M}_{ \Gamma_{I\infty} }})$. Here we follow  the
notation and results of \cite{GR}. For the definition and study of
cusps forms, readers are referred to the works of Gekeler, Goss or
the Habilatationshrift of Gebhard B\"ockle.

From the above lemma we obtain an additive version of
Stickelberger's theorem for $n=2$:

{\bf {Theorem 2}} If the   cardinal of the group $Pic ({\mathcal
E}({I\infty}))$ is $\infty$, there exists a cusp form for
$\Gamma_{I\infty}$ that is annihilated by $\tilde T (d\cdot 2)+
\tilde T (d\cdot 2-1)+\cdots + \tilde T (1)+ Id $.

$\tilde T (j)$ denotes the linear operator given by j-hecke acting
on the cusp forms.

From theorem $1$, we  also  obtain
    ideal class group annihilators for the cyclotomic function fields in the spirit of
Stickelberger's theorem. We prove that the correspondence

$$  \sum_{i=0}^{nd} [\Gamma (Frob^{nd-i}) * (\sum_{ \underset{(I,q(t))=1 , deg(q(t))=i} {q(t) \text{ monic}\in
{\Bbb F}_q[t]  }}    \varphi(q(t),n)\cdot \Gamma(q(t)  )]$$ is
trivial on $Spec(K^\infty_I \otimes K^\infty_I)$. $\Gamma(q(t) )$
denotes the graphic of the element of $G_I$ associated with the
class of $q(t)$ in $ ({\Bbb F}_q[t]/I)^\times$. $\varphi(q(t),n)$
indicates the   number of submodules $N\subseteq {\Bbb
F}_q[t]^{\oplus n}$ such that:
$$ {\Bbb F}_q[t]^{\oplus n}/N\simeq {\Bbb F}_q[t]/q_1(t)\oplus \cdots \oplus {\Bbb F}_q[t]/ q_n(t)$$
with the product of the invariant factors $q_1(t)  \cdots  q_n(t)$
equal to $q(t)$.  This latter result can also obtained in a more
direct way by using the Euler  product of section 2.4 and this
result for $n=1$. Bearing in mind
  the analogy between Drinfeld varieties in positive
characteristic  and modular curves for number fields, I believe
that the interest of this work is the possible translation of
results to modular curves.

\bigskip
\bigskip
 {\bf List of notations }
\bigskip

${{\Bbb F}_q}  $ is a finite field with $q$-elements, ($q=p^m$).

\bigskip
$\otimes$   denotes $   \otimes_{{\Bbb F}_q} $.

 \bigskip
 $\o_{{\Bbb P}^1}$ denotes the ring sheaf of the scheme ${\Bbb P}^1$.

\bigskip
$R $ is an ${{\Bbb F}_q}  $-algebra.

\bigskip
$R^\times$ denotes the group of units   in a ring $R$.
\bigskip
$F$, denotes the Frobenius morphism.
\bigskip
If $s\in Spec(R)$, then we denote by $k(s)$ the residual field
associate with $s$.

\bigskip
${\Bbb P}^1_R$ denotes ${\Bbb P}^1\otimes R$.

\bigskip
Let $S$ be a finite subset of geometric points of ${\Bbb P}^1$. We
denote by ${\Bbb A}^S$   the adele group  outside $S$, and $O^S$
denotes the adeles within ${\Bbb A}^S$ without poles.

\bigskip
Let $M$ be a vector bundle over $ {\Bbb P}^1 $; $M(k)$ denotes $M
\otimes_{\o_{{\Bbb P}^1}} {\o_{{\Bbb P}^1} }(k\infty)$, $k\in
{\Bbb Z}$.

\bigskip
If $f:X\to X$ is a morphism of separated schemes, then $\Gamma
(f)$ denotes the graphic of $f $; $\Gamma (f)=\{(f(x),x); \quad
x\in X\}$.

\bigskip
$\vert {\Bbb P}^1\vert$ and $\vert I\vert$ denote the geometric
points of ${\Bbb P}^1$ and $Spec({\Bbb F}_q[t]/I)$, respectively.
\section{Elliptic sheaves and Hecke correspondences. Euler products }
 In this section, except for proposition
\ref{inf},  all results are valid  for any smooth, geometrically
irreducible  curve over ${\Bbb F}_q$  provided with a rational
point $\infty$, although we only consider the projective line
curve.

\subsection{ Elliptic sheaves }

In this section we   recall the definition of   elliptic sheaves
and  level structures  over an ideal $I\subset
 {\Bbb F}_q[t] $,
   \cite{ BlSt}, \cite{ Dr2}, \cite{LRSt}, \cite{Mu}.

\begin{defn}\label{drin} An elliptic sheaf of rank $n$ over $R$,  $E:=(E_{j},i_{j},\tau)$, is a commutative
diagram of vector bundles of rank $n$ over ${\Bbb P}^1_R$, and
injective morphisms of modules $\{i_h\}_{ h\in {\Bbb N}}$, $\tau$:

$$ \xymatrix @C=35pt{
E_1 \ar[r]^{i_1} & E_2 \ar[r]^{i_2} &
 {\cdots} \ar[r]^{i_{(n-1)}}&
E_n  \ar[r]^{i_n} &
\cdots \\
   E^\sigma_0\ar[u]^\tau \ar[r]^(.5){  i^\sigma_0} &
  E^\sigma_1 \ar[u]^\tau \ar[r]^(.6){  i^\sigma_1} &
{\cdots }\ar[u]^\tau \ar[r]^(.6){ i^\sigma_n} &
  E^\sigma_n \ar[u]^\tau \ar[r]& {\cdots}}$$

(here, $ E^{\sigma }_i $ denotes $(Id\times F )^*  E_i $),
satisfying:

a) For any $s\in Spec(R)$, we fix $deg((E_{j})_s)= j$.
 \bigskip \bigskip

b) For all $j\in {\Bbb  Z}$, $E_{j+n}=E_{j}(1)$. We can assume
that the $i_k$ are natural inclusions.
\bigskip \bigskip

c) $E_{j}+\tau(  E^\sigma_j )=E_{j+1}$ for all $j$.
\bigskip \bigskip

d) $\alpha^*(E_i/E_{i-1})$ is a   rank-one free module over $R$,
$\alpha$ being  the natural inclusion $\infty \times Spec(R)
\hookrightarrow   {\Bbb P}^1_R$.

 \end{defn}
\begin{rem}\label{trivial} From these properties, it may be  deduced that $h^0(E_{j})=n+j$ and $h^1(E_{j})=0$,
$j\geq -n$, c.f. \cite{ BlSt}, {\cite{Dr2}}.

Moreover, it is seen  that for the $R$-module $H^0( {\Bbb P}^1_R ,
E_j)$, ($j>-n$), there exists
    a basis $\{s,\tau  s,\cdots, \tau^{ n+j-1}    s\}$     with
$\tau   s:=\tau((Id\times F)^* s)$ and $\tau^h  s:=\tau((Id\times
F)^*\tau^{h-1}\cdot s)$.
\end{rem}

\begin{defn}
An    $I$-level structure, $\iota_I$, for the elliptic sheaf
$(E_j,i_j,\tau)$   is an $I$-level structure, $\iota_{j , I} $,
for each vector bundle $ E_j$ compatible with the morphisms
$\{i_j,\tau\}$. i.e., $ (\iota_{j+1, I}) .i_j= \iota_{j  , I} $
and $(\iota_{j+1, I})\cdot \tau= (Id\times F)^*(\iota_{j , I}) $.
We denote by $(E,\iota_I)$  an elliptic sheaf with an $I$-level
structure.
 \end{defn}
Recall that an $I$-level structure for a vector bundle $E_j$ over
${\Bbb P}^1_R$  is a surjective morphism of modules $E_j\to
(\beta_*(R[t]/I))^{\oplus n} $, where
$\beta:Spec(R[t]/I)\hookrightarrow  {\Bbb P}^1_R $  is the natural
inclusion.

The elliptic sheaf $(E_j,i_j,\tau)$ defined over $R$  gives a
$\tau$-sheaf, $R\{\tau\}=\oplus_{i=0}^{\infty}R\cdot \tau^i$,
($\tau\cdot b=b^{q}\cdot \tau$). One can identify:
$$H^0( {\Bbb P}^1_R  ,E_j)=\oplus_{i=0}^{n+j-1}R\cdot
\tau^i  s,$$
  and in this way  $R\{\tau\}$ is isomorphic to the graded $
  R[t]$-module:
$$ \cup_{i=0}^{\infty}H^0( {\Bbb P}^1_R  ,E_j(i)).$$

\begin{rem}\label{dt}
By taking the determinant of $(E,\iota_I)$ we obtain  an elliptic
sheaf of rank $1$, $( det( E_{j}),   det( i_{j}),   det( \tau))$,
with an $I$-level structure $ det( \iota_I) $. This determinant is
studied in detail in {\cite {Ge}}.

The $\tau$-sheaf associated with $( det( E_{j}),   det(i_{j}),
\tau_{det})$ is
$$R\{ \tau_{det}\}:=\oplus_{i=0}^{\infty}R \cdot
\tau_{det}^i,$$ with $\tau_{det}^i=\tau^{i }\wedge \tau^{i+1 }
\wedge\cdots \wedge \tau^{n+i -1}$, ($i\geq 0$). Moreover,
$\wedge^n R\{\tau\}=R\{\tau_{det}\}$ as $R[t]$-modules.
\end{rem}

We denote   $ det( E_{j}) $ by $  L_j$; so $deg(L_j)=j$; recall
that $deg(E_j)=j$.

\begin{propo}\label{inf} With the above notations, if $r_n-r_1\geq n$
then:
$$\tau^{r_1  }\wedge \tau^{r_2 } \wedge\cdots \wedge \tau^{r_n }\in H^0({\Bbb P}^1_R,L_{r_n-n}) .$$
\end{propo}
\begin{proof} Since
$$\tau_{det}^{r_1}(1\wedge \tau^{r_2-r_1 } \wedge\cdots \wedge \tau^{r_n-r_1  })=\tau^{r_1  }\wedge \tau^{r_2 }
\wedge\cdots \wedge \tau^{r_n }\in H^0({\Bbb
P}^1_R,L_{r_n-n}),\quad  (0\leq r_1\leq \cdots\leq r_n)  $$ it
suffices to prove the result for $r_1=0$.

We proceed by induction over $r_n$. For $r_n=n$, we have to prove
that
$$1\wedge \tau^{r_2  } \wedge\cdots \wedge \tau^{n }\in H^0({\Bbb P}^1_R,L_{ 0}),$$

  $t\cdot a_n=a_0+a_1\cdot \tau+\cdots +  \tau^n$  for   $a_i\in R$, and therefore there exists $c\in R$ with
$$1\wedge \tau^{r_2  } \wedge\cdots \wedge \tau^{n }=c\cdot(1\wedge \tau^{ 2  } \wedge\cdots \wedge \tau^{n-1 }).$$
Recall that $(0\leq r_2\leq \cdots\leq  n)$ and we conclude since
$$1\wedge \tau^{ 2  } \wedge\cdots \wedge \tau^{n-1 }\in H^0({\Bbb P}^1_R,L_{ 0}).$$

Let us now assume the theorem is true for $k<r_n$. Let take us
$r_n=l+n$. Thus,
$$1\wedge \tau^{r_2  } \wedge\cdots \wedge \tau^{r_n }=1  \wedge\cdots \wedge \tau^{r_{n-1} }\wedge
 (t\cdot a_n\tau^{r_n-n  }-a_0\tau^{r_n-n  }+a_1\cdot \tau^{r_n-n+1  }+\cdots +  a_{n-1}\tau^{r_n -1  }).$$
Since
$$1  \wedge\cdots \wedge \tau^{r_{n-1} }\wedge\tau^{r_n-i  }\in  H^0({\Bbb P}^1_R,L_{ r_n -n }), \text { with }i\geq 1,$$
because
$$1 ,\cdots, \tau^{r_{n-1}
},\tau^{r_n-i  }\in  H^0({\Bbb P}^1_R,E_{ r_n -n }),$$

it suffices to prove that
$$(t\cdot a_n-a_0)\cdot (1  \wedge\cdots \wedge \tau^{r_{n-1} }\wedge
 \tau^{r_n-n  })\in H^0({\Bbb P}^1_R,L_{ r_n -n }).$$

Let consider us $k:=max\{r_n-n, r_{n-1}\}$. If $n\leq k$ then we
finish by induction, because $k+1 \leq r_n  $. When $k\leq n-1$,
it suffices to prove that
$$(t\cdot a_n-a_0)\cdot (1  \wedge\cdots \wedge \tau^{ {n-2} }\wedge
 \tau^{n-1  })\in H^0({\Bbb P}^1_R,L_{ r_n-n }).$$
This is true because we are in the case    $ r_n-n\geq 1$.

\end{proof}

\bigskip

\subsection{ $\infty$-Level structures }

We shall now define level structures at  $\infty \in {\Bbb P}^1  $
over elliptic sheaves of rank $1$. To do so, we   take  into
account the results of \cite{An1} 6.1.1.
 We take
$t^{-1}$ as a local uniformizer at $\infty$.

  The composition   of the epimorphism
$$\o_{{\Bbb P}^1}(k )\to\o_{{\Bbb P}^1}(k )/\o_{{\Bbb P}^1}( k-1  )$$
with the isomorphism induced by the multiplication by $ t^{1-k}$
$$\o_{{\Bbb P}^1}(k )/\o_{{\Bbb P}^1} ( k-1 )\simeq \o_{{\Bbb P}^1} /\o_{{\Bbb P}^1}( -1),$$
gives us an $\infty    $-level structure  over  $\o_{{\Bbb P}^1}(k
)$.

\begin{defn}\label{infinity} An $\infty$-level structure for a rank-$1$ elliptic sheaf, $(  L_j,i_j,\tau   )$, over
$R$  is an $\infty$-level structure $ (   L_0, \iota_\infty)$ such
that the diagram

$$(*) \xymatrix {    L^\sigma_0    \ar[dr]^{       \iota^\sigma_\infty   }    \ar[r]^{\tau } &
 L_0(1)\ar[d]^{t^{-1} \cdot \iota_\infty     }
\\ &  \gamma_*(R[t^{-1}]/t^{-1} R[t^{-1}])} $$
is commutative. Here,   $\gamma:Spec(R[t^{-1}]/t^{-1}
R[t^{-1}])\hookrightarrow  {\Bbb P}^1_R $, is the natural
inclusion.
\end{defn}

We denote by, ${\mathcal E}_n^{I  }$ and ${\mathcal E}_n^{I
\infty}$ the moduli of elliptic sheaves with $I$-level structures,
$(E,\iota_I)$, and with $I+\infty$-level structures, respectively.
Where to give an $\infty$-level structure for $E$, $\iota_\infty$,
is to give an $\infty$-level structure for  the rank-$1$-elliptic
sheaf $ det( E) $. Henceforth, we denote with $(E,\iota_{I
\infty})$ an element
 $(E,\iota_I,\iota_\infty)\in {\mathcal E}_n^{I \infty}$.   There exists a morphism, $z:{\mathcal E}_n^{I \infty}\to
 Spec({\Bbb F}_q[t])  $, called the zero morphism, that is defined
 by:
$$z  (E,\iota_{I \infty})=supp(E_0/\tau(E_{-1}^\sigma))=supp(det(E_0)/\tau_{det}(det (E_{-1})^\sigma).$$

Bearing in mind the antiequivalence between elliptic sheaves and
Drinfeld modules, one can construct a ring
  $ {\mathcal B}_n^{I  }$ of dimension $n$ such that $Spec( {\mathcal B}_n^{I  })= {\mathcal E}_n^{I  }$. For
these results see
   \cite{ Dr1},
\cite{ Dr2}, \cite{Lm}, \cite{Mu}.

For $n=1$, it is not hard to obtain a ring  $ {\mathcal B}_1^{I
\infty }$, such that
$$ Spec( {\mathcal B}_1^{I \infty })={\mathcal E}_1^{I \infty}.$$
Moreover, the morphism of forgetting the $\infty$-level structure
$${\mathcal E}_1^{I \infty}\to {\mathcal E}_1^{I  }$$
is \'{e}tale, outside $\vert I\vert $. In the following remark we
calculate   ${\mathcal B}_1^{I \infty }$, explicitly.

\begin{rem} We consider the  rank-$1$-Drinfeld module
$\phi_t=a  \tau +\bar t$, defined over ${\Bbb F}_q[\bar t,a]$. We
shall now study what   an $\infty$-level structure for the
Drinfeld module $\phi $ is.

Let us consider a rank-$1$ elliptic sheaf, $(  L_j,i_j,\tau   )$,
associated with $\phi $ and let $\iota_\infty$ be an
$\infty$-level structure for $(  L_j,i_j,\tau   )$. We have the
morphisms of modules:
$$\iota_\infty:L_0\to  {\Bbb F}_q[\bar t,a] [t^{-1}]/t^{-1} {\Bbb F}_q[\bar t,a] .$$
We choose $s$ with  $H^0(L_0)=<s>$. Note that $s$ is a generator
of the line bundle $L_0$.
 We set $\iota_\infty(s)=\lambda$, and hence
$$\iota^{\sigma}_\infty:L_0^\sigma\to {\Bbb F}_q[a,\bar
t][t^{-1}]/t^{-1}{\Bbb F}_q[a,\bar t]$$ gives
$\iota^{\sigma}_\infty(s)=\lambda^q$. Also,
$$t^{-1}\cdot \iota
_\infty:L_0(1)\to {\Bbb F}_q[\bar t,a][t^{-1}]/t^{-1}{\Bbb
F}_q[a,\bar t]$$ is such that $t^{-1}\cdot \iota
_\infty(\tau(s))=t^{-1}\tau s=a^{-1}$, because $ t\cdot s =a \cdot
\tau  s  +\bar t\cdot s$ and  $\bar t\cdot t^{-1} =0$ as element
 in $${\Bbb F}_q[a,\bar t][t^{-1}]/t^{-1}{\Bbb F}_q[a,\bar
t].$$ Therefore, the above diagram is commutative  if and only if
$\lambda^q=\lambda\cdot a^{-1}$. Thus, we can choose $\bar s:=
\lambda\cdot s\in H^0(M_0)$ such that $\iota _\infty(\bar s)=1$.
Therefore $t\cdot \bar s=\tau s+\bar t$, and we obtain the
Drinfeld module $\bar \phi_t=  \tau +\bar t$ isomorphic to $
\phi_t$. It is not hard to see that ${\mathcal B}_1^\infty ={\Bbb
F}_q[ \bar t]$.

We set $I=p(t){\Bbb F}_q[  t]$, where $d+1$ is the degree of
$p(t)$.
$$ \bar\phi_{p(t)}=c_{d+1}\tau^{d+1}+\cdots+c_1\tau +  p(\bar t); \quad c_i\in {\Bbb F}_q[\bar t].$$

We have that ${\mathcal B}_1^{I \infty  }={\Bbb F}_q[\bar t,p(\bar
t)^{-1},\delta]$
 with $\delta$ an element of a closed field of ${\Bbb F}_q(t)$
 verifying:
  $$  \phi_{p(t)}(\delta)= \delta^{q^{d+1}}+\cdots+c_1\delta^{q } +r(\bar t)\delta=0,$$
and  $  \phi_{h(t)}(\delta)\neq  0$ with $h(t)$ a proper  divisor
of $  r(t)$. The $I$-level structure for $(  L_j,i_j,\tau   )$ is
given by
$$\iota_I(\bar s)= \bar\phi_{t^{r-1}}(\delta)
+\bar\phi_{t^{d}}(\delta)t+\cdots+\bar\phi_{t^{0}}(\delta)t^{d}\in
{\Bbb F}_q[\bar t,\delta][t]/p(  t).$$
\end{rem}

The   morphism  ${\Bbb F}_q[t]\hookrightarrow  {\mathcal B}_1^{I
\infty}$  ($t\to \bar t$)  gives us the Galois extension
$K_I/{\Bbb F}_q(t)$   of group  $({\Bbb F}_q[t]/I)^\times$.

By considering
$$det  (E,\iota_{I }):=(det(E),det(\iota_{I  })),$$
 and the determinant morphism
$$det:{\mathcal E}_n^{I } \to  {\mathcal E}_1^{I },$$
we obtain:
$${\mathcal E}_n^{I \infty}={\mathcal E}_n^{I  }\times_{{\mathcal E}_1^{I  }} {\mathcal E}_1^{I \infty},$$
and therefore, ${\mathcal E}_n^{I \infty}$ is an  affine scheme of
finite type over ${\Bbb F}_q$. It is smooth because the projection
$${\mathcal E}_n^{I  }\times_{{\mathcal E}_1^{I  }} {\mathcal E}_1^{I \infty}\to {\mathcal E}_n^{I  }$$
is \'{e}tale  since  ${\mathcal E}_1^{I \infty}\to {\mathcal E}_1^{I
}$  is also \'{e}tale . Note that ${\mathcal E}_1^{I \infty}$ is
defined over  ${\Bbb P}^1\setminus \vert I \vert \cup \infty$.


\subsection{ Hecke correspondences}
We consider,  $J_1\subseteq \cdots \subseteq J_n$, a chain of
ideals of ${\Bbb F}_q[t]$ coprime to $I$ and $S=\vert I \vert \cup
\{\infty\}$.

 Let $(E,\iota_{I \infty})$ be an elliptic sheaf defined over
$R$ with level structures on $I$ and on $\infty$ and with zero
outside $\vert J_1\vert$. We denote   by ${\mathcal E}_{n, \vert
J_1\vert}^{I \infty}$ the moduli scheme
$${\mathcal E}_n^{I \infty} \times_{{\Bbb P}^1}  ({{\Bbb P}^1}\setminus \vert J_1\vert) ,$$
where the fibred product is obtained from the zero morphism
$z:{\mathcal E}_n^{I \infty}\to {{\Bbb P}^1}$ and the natural
inclusion ${{\Bbb P}^1}\setminus \vert J_1\vert\ \hookrightarrow
{{\Bbb P}^1}$.

  We denote by
$$T(J_1, \cdots , J_n)\subset  {\mathcal E}_{n, \vert J_1\vert}^{I \infty} \times  {\mathcal E}_{n, \vert J_1\vert}^{I
\infty}$$  the Hecke correspondence, which is given by the pairs
$$[(E,\iota_{I \infty}),(\bar E,\bar\iota_{I \infty})]\in
 {\mathcal E}_{n, \vert J_1\vert}^{I \infty}\times {\mathcal E}_{n, \vert J_1\vert}^{I \infty},$$

$ E$ being  a sub-elliptic sheaf  of $ \bar E $ such that for each
$s\in Spec(R)$ we have
$$\bar E_s/  E_s\simeq   k(s)[t] /J_1  \oplus \cdots \oplus k(s)[t] /J_n .$$
The  $I+\infty$-level structure,   $ \iota_{I \infty} $, defined
over    $E$ is the
  composition $\bar\iota_{I\infty}\cdot\rho$,  $\rho$ being the inclusion $E\subset \bar E$.

We shall now describe the Hecke correspondences in an adelic way.
To do so,  consider   $(\bar E,\bar \iota_{I \infty})$ defined
over an algebraic closed field $K$.

We denote by
$$\pi_1, \pi_2:{\mathcal E}_{n, \vert J_1\vert}^{I \infty}\times {\mathcal E}_{n, \vert J_1\vert}^{I \infty}\to {\mathcal
E}_{n, \vert J_1\vert}^{I \infty}$$   the  natural projections.
There exists a bijection between  the sets:
$$\pi_1(\pi_2^{-1} (\bar E,\bar \iota_{I\infty})\cap T(J_1, \cdots , J_n) ), \quad \pi_2(\pi_1^{-1} (\bar E,\bar
\iota_{I\infty})\cap T(J_1, \cdots , J_n) ),$$ and
 the ${\Bbb F}_q[t]$-submodules $M$ and $\bar M$,
$$ M\subseteq {\Bbb F}_q[t]^{\oplus n}  \subseteq \bar M$$
with
$$ {\Bbb F}_q[t]^{\oplus n}/M\simeq \bar M/{\Bbb F}_q[t]^{\oplus n} \simeq {\Bbb F}_q[t] /J_1  \oplus \cdots \oplus {\Bbb
F}_q[t] /J_n
$$ respectively. These sets have the same cardinal, which we denote by $d(J_1, \cdots , J_n)$.

In the following proposition,  $ {\mathcal Cht}_{n, \vert
J_1\vert}^{I  }$, denotes the stack of shtuckas of rank $n$ with
zeroes outside  $\vert J_1\vert$ and level structures over $I$.
c.f. \cite{Lf}
\begin{propo}$T(J_1, \cdots , J_n)$ is a closed subscheme of $  {\mathcal E}_{n, \vert J_1\vert}^{I \infty}\times
{\mathcal E}_{n, \vert J_1\vert}^{I \infty}$. Moreover, the
morphisms  $ \pi_1, \pi_2$ restricted to $T(J_1, \cdots , J_n)$
are     \'{e}tale morphisms. We denote these morphisms by $ \bar
\pi_1, \bar\pi_2$, respectively.
\end{propo}
\begin{proof}  Let consider us the morphism:
$$\mathfrak e:{\mathcal E}_{n, \vert J_1\vert}^{I  }\to {\mathcal
Cht}_{n, \vert J_1\vert}^{I  },$$ defined by $$\mathfrak e
(E,\iota_{I }):=((\xymatrix {    E_{-1}    \ar[r]^i      &  E_0
 & E^\sigma_{-1} \ar[l]_\tau  }),\iota_{I }). $$
 c.f.
\cite{Dr3}, (pag.109).

The Hecke correspondences $\Gamma^n(g) $ defined in \cite{Lf}
(section I, 4) are closed substacks within $ {\mathcal Cht}_{n,
\vert J_1\vert}^{I  }\times {\mathcal Cht}_{n, \vert J_1\vert}^{I
}$. Let us consider $g\in Gl_n({\Bbb A}^{S})$ such that
$$\oplus^n O^{S}/g(\oplus^n O^{S})\simeq {\Bbb F}_q[t] /J_1  \oplus \cdots \oplus {\Bbb F}_q[t] /J_n$$
as modules. In this way,
$$ T^I(J_1, \cdots , J_n) :=(\mathfrak e \times \mathfrak e )^*\Gamma^n(g) $$
  is a closed subscheme of ${\mathcal E}_n^{I  }\times
{\mathcal E}_n^{I}$, where $T^I(J_1, \cdots , J_n)$ denotes the
Hecke correspondence
$$(\pi_\infty \times \pi_\infty)(T (J_1, \cdots , J_n))\subset {\mathcal E}_{n, \vert J_1\vert}^{I  }\times
{\mathcal E}_{n, \vert J_1\vert}^{I  },$$ $\pi_\infty:{\mathcal
E}_{n, \vert J_1\vert}^{I \infty}\to {\mathcal E}_{n, \vert
J_1\vert}^{I  }$ being  the morphism of forgetting the
$\infty$-level structure. Now, $T (J_1, \cdots , J_n) $ is the
closed subscheme given by   the pairs:
$$ [(E,\iota_{I \infty}),(\bar E,\bar\iota_{I \infty})]\in (\pi_\infty \times \pi_\infty)^{-1}T^I(J_1, \cdots , J_n),$$
such that
$$ \xymatrix {    det(E)   \ar[dr]^{       \iota_\infty   }    \ar@{^{(}->}[r]^{det(\rho)} &
det(\bar E)\ar[d]^{\bar \iota_\infty    }
\\ &   \gamma_*(R[t^{-1}]/t^{-1} R[t^{-1}]) } $$
is commutative. $det(\rho):det(E) \hookrightarrow det(\bar E) $ is
the determinant of the injective morphism  given between the
elliptic sheaves $\rho:E\hookrightarrow \bar E$.

Because
$$T^I(J_1, \cdots , J_n)=\Gamma^n(g)\times_{{\mathcal Cht}_{n, \vert J_1\vert}^{I}}{\mathcal E}_{n, \vert J_1\vert}^{I},$$
and since the projections $p_i:\Gamma^n(g)\to {{\mathcal Cht}_{n,
\vert J_1\vert}^{I}}$, ($i=1,2$) are \'{e}tale morphisms, we have that
the two    projections from $ T^I(J_1, \cdots , J_n)$  to
${\mathcal E}_{n, \vert J_1\vert}^{I}$ are \'{e}tale morphisms.

We conclude that $ \bar \pi_1, \bar\pi_2$ are \'{e}tale morphisms
because
$$T (J_1, \cdots , J_n)=T^I(J_1, \cdots , J_n)\times_{{\mathcal E}_{n, \vert J_1\vert}^{I }}{\mathcal E}_{n, \vert
J_1\vert}^{I\infty}.$$ They are morphisms  of   degree    $d(J_1,
\cdots , J_n)$.

\end{proof}

\bigskip
The formal sum of Hecke correspondences gives a commutative ring
where the product is the composition of correspondences. This ring
is isomorphic to the commutative ring
$$C_c(K\backslash Gl_n({\Bbb A}^{S  })/K)$$
 of  ${\Bbb Z}$-valued
continuous  functions over $Gl_n({\Bbb A}^{S})$, invariant by the
action of $K:=Gl_n(   O^{S })$ on the left and on the right over
$Gl_n({\Bbb A}^{S })$ and with compact support. The product is the
convolution product. This isomorphism  sends the correspondence
$T(J_1, \cdots , J_n)$  to the characteristic function over the
open compact  subset:
$$Gl_n(   O^{S})\cdot  (\mu_{J_1},\cdots,\mu_{J_n})\cdot Gl_n(   O^{S }),$$
with $\mu_{J_i}\in {\Bbb A}^S$,
  given by the element $q_i(t)$ with   $J_i=q_i(t){\Bbb F}_q[t]$, and $(\mu_{J_1},\cdots,\mu_{J_n})$ the diagonal
matrix in $Gl_n({\Bbb A}^{S })$, where the diagonal is given by
$(\mu_{J_1},\cdots,\mu_{J_n})$.

We denote by $T(m)$ the   correspondence defined by the formal sum
of the Hecke correspondences $T(J_1, \cdots , J_n)$, where $
\sum^n_{i=0}dim_{{\Bbb F}_q}{\Bbb F}_q[t]/J_i=m$.
\bigskip

As in the number field case, one can consider Hecke
correspondences as operators over the abelian group of formal sums
of ${\Bbb F}_q[t]$-submodules, $N$, of rank $n$ of ${\Bbb
F}_q(t)^{ \oplus n  }$ (=lattices of ${\Bbb F}_q(t)^{ \oplus n
}$). One defines:
$$T(J_1,\cdots, J_n)(N)=\sum \bar N,$$
where $\bar N$ runs over the submodules of $N$, satisfying:
$$N/\bar N \simeq  {\Bbb F}_q[t]/ J_1   \oplus \cdots
\oplus {\Bbb F}_q[t]/ J_n.$$ In this way $T(m)(N)=  \sum_{\bar
N\subset    N}  \bar N$, where $dim_{{\Bbb F}_q}   N/\bar N=m$.

A more rigorous presentation of this section can be found  in
\cite{Lf} and \cite{Lm}.


\subsection{ Euler products}

A  generalization
 for the non-abelian case of the $S $-incomplete $L$-function
evaluator at $s=0$, (cf: \cite{H1}, \cite{Ta}) is studied in
\cite{H2}. In this section we address the issue in another way.

${\mathcal E}_{n,\vert {{\Bbb P}^1}\vert }^{I \infty}$ denotes the
moduli scheme of elliptic sheaves of rank $n$ with level
structures over $I$ and $\infty$ and with zero outside    $\vert
{{\Bbb P}^1}\vert$.

Let   $x\in \vert {{\Bbb P}^1}\vert\setminus S$, ($S= \vert I\vert
\cup \{ \infty\}$) and let $t_x$ be a local uniformizer for $x$.
We consider the diagonal  matrix
$$(\mu_x,\overset{\overset j \smile} \cdots,\mu_x,1,\cdots,1)\in Gl_n({\Bbb A}^{S  }),$$
  $\mu_x$ being the adele within ${\Bbb A}^{S  }$ such that it is $1$ over each place of $\vert {{\Bbb P}^1} \vert \setminus
S\cup \{x\}$ and $t_x$ over $x$. We denote by $\sigma^x_j$, $1\leq
j \leq n$, the Hecke correspondence over ${\mathcal E}_{n,\vert
{{\Bbb P}^1}\vert  }^{I \infty}$  given by  the characteristic
function of
$$Gl_n(   O^{S})\cdot  (\mu_x,\overset{\overset j \smile} \cdots,\mu_x,1,\cdots,1)\cdot Gl_n(   O^{S}).$$

In the following lemma, for easy notation we assume that
$deg(x)=1$ and $t_x$ is a local parameter for $x$. ${\mathfrak
m}_x$ is the maximal ideal associated with $x$.

One can find a proof of the next lemma in, \cite{Sh} Th 3.21. More
or less, we   repeat that proof.
\begin{lem}\label{euler} We have:
$$\frac{1}{ 1-\sigma^x_1\cdot z +  q \sigma^x_2\cdot z^{2 }
- q^{3 }\sigma^x_3\cdot z^{3 }+\cdots +              (-1)^n
q^{n(n-1)/2 } \sigma^x_n\cdot z^{n }}
   =\sum_{m \geq 0} T^x(m)\cdot z^m,   $$
where
$$T^x(m)\subset  {\mathcal E}_{n,\vert {{\Bbb P}^1}\vert }^{I \infty}\times  {\mathcal E}_{n,\vert {{\Bbb P}^1}\vert }^{I
\infty}$$ denotes   sum of the Hecke correspondences
 $T({\mathfrak m}_x^{r_1}, \cdots , {\mathfrak m}_x^{r_n})$, ($r_1\geq\cdots \geq r_n\geq 0$), and    $r_1+\cdots +r_n=m$.

\end{lem}

\begin{proof} We shall model this proof as in \cite{Ln}. It suffices to prove that for each $r\in {\Bbb N}$ we have  "Newton's" formulas
$$P:=T^x(r) -   T^x(r-1) \cdot \sigma^x_1 +  q      \sigma^x_2 \cdot T^x(r-2)
- \cdots +              (-1)^n  q^{n(n-1)/2 }
T^x(r-n)\cdot\sigma^x_n =0$$ by denoting $T^x(0)=1$ and $T^x(l)=0$
for $l<0$.

 To accomplish this, we
consider Hecke correspondences as operators over the formal
abelian group of lattices,  $\bar N$ and $N$ being   lattices
 with $\bar N\subseteq N$, $dim_{{\Bbb F}_q} N/\bar N=r$ and $N/\bar N$   concentrated over $x$. We shall prove that  the
multiplicity of $\bar N$ in the formal sum $P(N)$ is $0$.

$\sigma^x_j(N)=\sum N'$, where $N'$ belong to the sublattices
$N'\subset N$ with
$$N/N'\simeq {\Bbb F}_q[t]/{\mathfrak m}_x  \oplus \overset {\overset j \smile} \cdots \oplus
{\Bbb F}_q[t]/{\mathfrak m}_x$$ or, equivalently, the vector
subspaces, $N'$ of codimension $j$, of  $N/ {\mathfrak m}_x\cdot
N$.

If we denote $h:=dim_{{\Bbb F}_q} N/({\mathfrak m}_x\cdot N+\bar
N)   $, then the number of lattices, "$N'$" such that $\bar
N\subset N'$, is given by  the number of ${\Bbb F}_q$-subvector
spaces in $N/({\mathfrak m}_x\cdot N+\bar N)$ of codimension $j$.
This number is given by the $q$-combinatorial number
$$\left(\begin{matrix}     h
 \\ j  \end{matrix}\right)_q=:\frac{(q^h-1)\cdots (q^{h-j+1}-1)}{ (q^j-1)\cdots (q-1)} $$
for $j\leq h$, and  $\left(\begin{matrix}     h
 \\ j  \end{matrix}\right)_q=:0$ for either $ j<0 $ or  $ j>h$.

We conclude the lemma bearing in mind the relation
$$\left(\begin{matrix}     h
 \\ h  \end{matrix}\right)_q -   \left(\begin{matrix}     h
 \\ h-1  \end{matrix}\right)_q  +  q  \cdot \left(\begin{matrix}     h
 \\ h-2  \end{matrix}\right)_q
- \cdots +              (-1)^n  q^{n(n-1)/2 }   \cdot
\left(\begin{matrix}     h
 \\ h-n  \end{matrix}\right)_q=0.$$
I have taken  this formula from Appendix D, \cite{Lm}, (cf.
\cite{Ma}).

\end{proof}

\begin{thm}
If we denote
$$L^x:= \frac{1}{ 1-\sigma^x_1\cdot z +  q \sigma^x_2\cdot z^{2 }
- q^{3 }\sigma^x_3\cdot z^{3 }+\cdots +              (-1)^n
q^{n(n-1)/2} \sigma^x_n\cdot z^{n }},  $$ then
$$\prod_{x\in \vert {\Bbb P}^1 \vert\setminus S} L^x=\sum_{m \geq 0} T (m)\cdot z^m .$$
\end{thm}
\begin{proof} It suffices to bear in mind the above lemma and that if $J_1$ and $\bar J_1$ are
ideals coprime  within $A$,
 then:
$$T(J_1,\cdots, J_n)\cdot T(\bar J_1,\cdots, \bar J_n)=T(J_1\cdot \bar J_1,\cdots, J_n\cdot \bar J_n).$$
\end{proof}


\section{ Isogenies and Hecke correspondences}
  Here,   we study the isogenies between Drinfeld modules(=elliptic sheaves),
  \cite{Gr2},
to establish the relation between the above Euler   products and
isogenies between elliptic sheaves.
\subsection{ Isogenies for elliptic sheaves}

\begin{defn} An isogeny, ${  \Phi}$, of degree $m\in {\Bbb N}$ between two elliptic sheaves with $I
$-level structures $(E,\iota_{I \infty}),(\bar E,\bar\iota_{I
\infty} )$ and $\infty $-level structures for $det( E)$ and
 $det(  \bar E)$    is a morphism of modules
$\Phi_i:E_i \to \bar E_{i+m}$, for each $i$, with   $Im
(\Phi_i)\not\subset  \bar E_{i+m-1}$,  preserving the diagrams
that define the elliptic sheaves and their  level structures.
\end{defn}
If $E$ and $\bar E$ are defined over $R$, then to give an isogeny,
${  \Phi}:E\to \bar E$,    of degree $m$ is equivalent to giving a
morphism of $\tau$-sheaves $\phi:R\{\tau\}\to R\{\bar\tau\}$, such
that if $r(\tau)$ is a monic polynomial  with
$deg_\tau(r(\tau))=j$ then $deg_{\bar \tau} \phi(r(\tau)) = m+j$.

\begin{lem} Let $M$ and $N$ be vector bundles of rank $n$ over ${\Bbb P}^1_R$, and with $x$ a rational point of
${\Bbb P}^1$. If $f:M\to N$ is a morphism of modules such that its
restriction to $k(x)$
$$f_{\vert k(x)\otimes R}:M_{\vert
k(x)\otimes R}\to N_{\vert k(x)\otimes R}$$
 is an isomorphism, then $f$ is injective.
\end{lem}
\begin{proof} Let assume us that     $x$ is the rational point  $0\in {\Bbb P}^1$. We have the exact
sequence:
$$0\to K\to M\overset f \to N.$$
If we prove that $K_{(\vert {\Bbb P}_1\setminus \{\infty\})\otimes
R}=0$ then we  conclude. Let
$$0\to \hat K\to \hat M\overset {\hat f} \to \hat N$$
be
  the completion of the above exact sequence along the ideal $t  R[t] $. By
hypothesis, $f_{\vert k(x)\otimes R}$ is an isomorphism. One
deduces that ${\hat f}$ is also an isomorphism   and hence $\hat
K=0$, since
$$Spec_{\text{maximal}}(R[[t]])=0\times  Spec_{\text{maximal}}(R),$$
  and in view of the
Nakayama lemma. If we prove that the natural morphism $g:K\to \hat
K$ is injective we conclude. By the Krull Theorem,  if $g(a)=0$
then there exists $1+t\cdot q(t) \in R[t]$ such that $(1+t\cdot
q(t) )\cdot a=0$.  However, the homothety morphism given by $
1+t\cdot q(t) $ over $R[t]$ is injective and therefore it is also
injective over $M$ because $M$ is locally free. Hence,   $a=0$.

\end{proof}

\begin{lem}\label{un} Assuming the above notations, if $\Phi$ is an isogeny   of degree $m\leq n\cdot d  $
between $(E,\iota_{I \infty} ),(\bar E,\bar\iota_{I \infty} )$
then $\Phi$ is injective and it is the only isogeny between these
elliptic sheaves with level structures. Moreover, there exists $
r\in {\Bbb N}$, maximum ($r\leq n\cdot d $), such that $\phi
(R\{\tau\} )\subseteq  R\{\bar\tau\}\cdot \bar \tau^r $.
\end{lem}
\begin{proof}
We can assume that the elliptic sheaves are defined over an
 ${\Bbb F}_q$-algebra   $R$. In this way, the injectivity is deduced from the
above lemma. We denote by $I$ indistinctly the ideal within ${\Bbb
F}_q[t]$ as the ideal sheaf  within $\o_{{\Bbb P}^1}$

Let $\Phi'$ be another isogeny; $\Phi-\Phi'$ defines a morphism
$E_0\to  I\cdot \bar E_{ n\cdot d}  $. Since   $E_0$ is generated
by its   global sections, and since $deg(I\cdot \bar E_{ n\cdot d}
))=-n$, because $deg(I)=d+1$, we have $h^0( I\cdot \bar E_{ n\cdot
d} )= 0$. In this way we have that $\Phi=\Phi'$.

The last assertion of the lemma is evident.

\end{proof}

We consider $\vert {\Bbb P}_1\vert_{nd}$, the subset of geometric
points of ${\Bbb P}_1$, of degree less than or equal to $nd$.
 ${\mathcal E}_{n,\vert {{\Bbb P}^1}\vert_{nd} }^{I \infty}$ denotes the moduli scheme of elliptic sheaves of rank $n$ with
level structures over $I$ and $\infty$ and with zero outside
$\vert {{\Bbb P}^1}\vert_{nd}$.

\begin{lem} With the above notations,  the  set
$$[(E,\iota_{I \infty} ),(\bar E,\bar\iota_{I \infty} )]\in
{\mathcal E}_{n,\vert {\Bbb P}_1\vert_{nd}}^{I \infty}\times
{\mathcal E}_{n,\vert {\Bbb P}_1\vert_{nd}}^{I \infty}, $$ such
that there exists an isogeny of degree $m\leq n\cdot d  $  between
$(E,\iota_{I \infty} )$ and $(\bar E,\bar\iota_{I \infty} )$ with
$r=0$, is given by the correspondence $  T(m) \subset {\mathcal
E}_{n,\vert {\Bbb P}_1\vert_{nd}}^{I \infty}\times {\mathcal
E}_{n,\vert {\Bbb P}_1\vert_{nd}}^{I \infty}. $
\end{lem}
\begin{proof}It is clear that a pair within $T(m)$ defines an isogeny of degree $m$ with the  required properties.
Moreover, the very lemma   asserts that  there only exists  one
isogeny of degree $m\leq nd$ between   two elliptic sheaves with
$I$-level structures. With this result, one deduces that if
$(E,\iota_{I \infty} )$ and $(E',\iota'_{I \infty} )$ are
subelliptic sheaves with level structures of $(\bar E,\bar\iota_{I
\infty} )$, by two different isogenies of degree $m$, then
 $(E,\iota_{I \infty} )$ is not isomorphic to $ (E',\iota'_{I \infty} )$.

On the other hand, if $\Phi:(E,\iota_{I \infty} )\to (\bar
E,\bar\iota_{I \infty} )$ is an isogeny with $r=0$ and degree $m$,
then by the serpent lemma  we have isomorphisms $ (Id\times
F)^*(\bar E_{i+m}/\Phi_i(E_i))\simeq \bar E_{i+m}/\Phi_i(E_i)$,
for each integer, $i$. Since the zeroes of the elliptic sheaves
considered
 are of degree $> nd$, we have
$$\bar E_{i+m}/\Phi_i(E_i)\simeq R[t] /J_1  \oplus \cdots \oplus R[t] /J_n,$$
where $ J_1 \subseteq  \cdots   \subseteq  J_n $ are ideals within
${\Bbb F}_q[t] $ coprime  to $I$ with $ \sum^n_{i=0}dim_{{\Bbb
F}_q}A/J_i=m$. Here, we have assumed that $(E,\iota_{I \infty} )$
and $ (\bar E,\bar\iota_{I \infty} )$ are defined over $R$.

\end{proof}

\begin{cor}\label{tri} The subset of pairs $(e,\bar e)\in {\mathcal E}_{n,\vert {\Bbb P}_1\vert_{nd}}^{I \infty}\times
{\mathcal E}_{n,\vert {\Bbb P}_1\vert_{nd}}^{I \infty}$ such that
there exists an isogeny of degree $n\cdot d$ is given by the
correspondence:
$$T(n\cdot d)+\Gamma(Fr)* T(n\cdot d-1)+\cdots +\Gamma(Fr^{n\cdot d-1})* T (1)+\Gamma(Fr^{n\cdot d }).$$
Here, $\Gamma(Fr^i)$ is given by the graph   of the
$q^i$-Frobenius morphism. $*$ denotes the product of
correspondences.

\end{cor}
\begin{proof}The elliptic sheaf associated with the $\tau$-sheaf, $R\{\bar \tau\}\cdot
\bar\tau^r$, is
$$[(Id\times F)^r]^*\bar E.$$
In view of the two last lemmas, the corollary is deduced bearing
in mind that between $E_0$ and $  [(Id\times F)^{nd+j}]^*\bar
E_{-j}$ there is no injective morphism for $j>0$, because
$deg(E_0)=0$ and $deg[((Id\times F)^{nd+j})^*\bar E_{-j}]= -j$.
\end{proof}

\subsection{ Trivial correspondences}

In this section we shall prove that the correspondence  of the
above corollary \ref{tri} is trivial.

\begin{propo} \label{sequence}Let $M$ be a   vector bundle over ${\Bbb P}^1_R $  of rank $n$
and degree $0$ where $h^0(  M(-1) )= h^1(  M(-1) )=0$, and with an
$    I$-level structure $\iota_{    I}$. Thus, we have that
$H^0({\Bbb P}^1_R , M )$ is a free  $R$-module of rank $n$, and
$M\simeq H^0({\Bbb P}^1_R , M )\otimes \o_{{\Bbb P}_1}$.

\end{propo}

\begin{proof} If   $x\in  Spec({\Bbb F}_q[t]/I) $ is a rational point, then   $h^0(  M(-x) )=
h^1(  M(-x) )=0$.
 Bearing in mind the morphism given by the $x$-level structure $   \iota_x:M\to (k(x)  \otimes R)^n$, we
obtain an isomorphism:
$$M/M(-x)\simeq  (R[t]/ {\mathfrak m}_x)^{\oplus n}.$$
Therefore, by taking global sections in the exact sequence of
$\o_{{\Bbb P}_1}\otimes R$-modules
$$0\to M(-x )\to M \to M/M(-x ) \to 0 $$
  we conclude.

The argument is valid when $Spec({\Bbb F}_q[t]/I)$ does not have
rational points because   $N$ is a $R$-free module if and only if
$N\otimes {\Bbb F}_{q^d}$ is a $R\otimes {\Bbb F}_{q^d}$-free
module.

  \end{proof}

If $(\bar M,\bar \iota_{    I})$ is an $I$-level structure then we
denote by $(\bar M(d),\bar \iota_{    I}(d))$ the $I$-level
structure over $\bar M(d)$ obtained from $\bar \iota_{    I}$ by
considering the natural inclusion $\bar M\subseteq \bar M(d)$.
Recall that $\infty \notin \vert I\vert$.

\begin{lem}\label{morfismo} If $(\bar M,\bar \iota_{    I}),(M,\iota_{    I}) $ are level structures
  over $R$, where $M$ and $\bar M$ satisfy the conditions of the above proposition, then there exists a   morphism of
vector bundles, $h: M\to \bar M(d ), $ whose diagram
$$ \xymatrix {    M   \ar[dr]^{  \iota_{   I}}    \ar[r]^{h} &
 \bar M(d )\ar[d]^{\bar \iota_{ I }(d)  }
\\ &(\beta_*(R[t]/I))^{\oplus n    } }$$
 is commutative. $h$ is said to be a morphism between  $(M,\iota_{    I}) $ and $(\bar M(d),\bar \iota_{    I}(d))$.
\end{lem}
\begin{proof} By choosing a base,
$\{s_1,\cdots,s_n\}$  for $H^0({\Bbb P}^1_R, M )$.

  $ H^0(\iota_{   I}):H^0(M)\to  (R[t]/I)^{\oplus n    }$   has the associated
  matrix:
$$ \Delta_0+\Delta_1t+\cdots +\Delta_{d }t^{d },  $$
  where $\Delta_i$  are $n\times n$-matrices with entries in $R$.

 We have  that
$$H^0({\Bbb P}^1_R, \bar M(d ) )=    \oplus^d_{i=0} H^0({\Bbb P}^1_R,
\bar M )\cdot t^i.$$ Moreover,
 bearing in mind that $deg( I)= {d+1}$ we also have that
$$H^0({\Bbb P}^1_R, \bar M(d) )\overset{H^0(\bar \iota_{ I})} \simeq (R[t]/I)^{\oplus n    },$$
because $h^0(I\cdot \bar M(d))=0$. Thus, we have that $H^0(h)$
must satisfy
$$  H^0(h)=A_0+   \cdots +A_{d }.t^{d }:=(H^0(\bar \iota_{ I }))^{-1}\cdot(\Delta_0+\Delta_1t+\cdots +\Delta_{d }t^{d }) ,
$$
where $A_i$ are $n\times n$-matrices with entries in $R$. We
conclude by bearing the above proposition in mind.
\end{proof}

The same arguments of lemma \ref{un} allow us to deduce that $h$
is unique.

By remark \ref{trivial},   if $E$ is an elliptic sheaf of rank
$n$, then   $E_0$ satisfies   the conditions of
  proposition
\ref{sequence}.

Let us consider the elliptic sheaves, defined over $R$, with $I
$-level structures $(E,\iota_{I \infty}),(\bar E,\bar\iota_{I
\infty} )$ and   $\infty $-level structures for $det( E)$ and
 $det(  \bar E)$.

\begin{lem}\label{condiciones} Let  $h $ be the morphism between vector bundles with level structures
$h :(E_0,\iota_{0,I}   )\to (\bar E_0(d) ,\bar\iota_{0,I}(d) )$
given in   lemma \ref{morfismo}, and let $ \iota_{\infty }$,
$\bar\iota_{\infty }$ be level structures at   $\infty$ over
$det(E)$ and $det(\bar E)$, respectively. Therefore,  there exists
an isogeny of degree $ n\cdot d$, $\Phi :(E,\iota_{I\infty} )\to
(\bar E,\bar\iota_{I\infty})$ with $\Phi_0=h $ if and only if
$det( h )$ is  a morphism for the level structures $(det(
E_0),\iota_{\infty}   )$,  $(det(  \bar E_0(d))
,\bar\iota_{\infty} )$, (i.e: $\bar\iota_{\infty}\cdot det(  h)=
\iota_{\infty}$),  and the morphism of $  R[t]$-modules given by
$h $, $h_A:R\{\tau\}\to R\{\bar\tau\}$ satisfies $deg_{\bar
\tau}(h(1))\leq n-2+nd ,\cdots ,deg_{\bar \tau}(h(\tau^{n-2}))\leq
n-2+nd$.
\end{lem}
\begin{proof}The direct way is trivial.

We prove the converse. Since $h$ is a morphism for the $I$-level
structures $(E_0,\iota_{0,I}   )$ and $ (\bar E_0(d)
,\bar\iota_{0,I}(d) )$, we have
  that
$$ h_A ( \tau)-\bar \tau \cdot  h_A(1)  ,\cdots ,h_A   (\tau^{n-1})-\bar \tau\cdot  h_A(
 \tau^{n-2}) \in I \cdot  R\{\bar \tau\}   .$$
Moreover, since   $deg(I)=d+1$, if $r(\bar \tau )\neq 0\in R\{\bar
\tau\}\cdot I$ then $deg_{\bar \tau}(r(\bar \tau ))\geq n(d+1)$.
Thus, by the hypothesis of the lemma we deduce that
$$ h_A ( \tau)-\bar \tau \cdot  h_A(1) =\cdots =h_A   (\tau^{n-1})-\bar \tau\cdot  h(
 \tau^{n-2}) =0.$$
Therefore, $h_A(\tau^i)=\bar \tau^j\cdot h_A(  \tau^k)$ for
$j+k=i$ and $1\leq i\leq n-1$.

Now, we prove the equalities
$$ (*) \quad deg_{\bar \tau}(h_A(  1))= nd  ,  deg_{\bar \tau}( h_A  (\tau   ))=nd+1   ,  \cdots ,  deg_{\bar \tau}(  h_A
(\tau^{n-1} ))=nd+n-1 .$$

We   consider the determinant elliptic sheaves $det(  E)$ and
$det( \bar E)$ and their  $\tau$-sheaves $R\{\tau_{det}\}$ and
$R\{ \bar \tau_{det})\}$, respectively. We have that
$$[det(  h_A )\cdot \tau_{det}- \bar \tau_{det} \cdot det( h_A )](  1\wedge\tau    \wedge
\cdots\wedge \tau^{n-2}     \wedge\tau^{n-1}      ) $$ is an
element of $R\{\tau_{det}\}$ of degree $\leq n d+1 $. However,  by
hypothesis
 $det(
h )$ is  a morphism for  $\infty$-level structures for elliptic
sheaves and therefore this element is of degree  $\leq nd$.

Because $h_A(\tau^i)=\bar \tau \cdot h_A(  \tau^{i-1})$, for
$1\leq i\leq n-1$,  the above element of $R\{\tau_{det}\}$ is
equal to
$$
 h_A  (\tau   )\wedge  \cdots\wedge  h_A  ( \tau^{n-1}     )\wedge
  [h_A \cdot \tau-\bar \tau\cdot h_A ](\tau^{n-1}  )  .$$

Since      $ deg_{\bar \tau}(h_A  (\tau^{n-1} ))\leq nd+n-1 $ and
$h_A(\tau^i)=\bar \tau^i\cdot h_A(1)$ for $1\leq i\leq n-1$, we
have the inequalities
$$ deg_{\bar \tau}(h_A(  1))\leq nd  ,  deg_{\bar \tau}( h_A  (\tau
))\leq nd+1   ,
 \cdots ,
 deg_{\bar
\tau}(  h_A  (\tau^{n-1} ))\leq nd+n-1. $$
 But
  $ \bar E_{0} (d)/h (E_0)$ is not concentrated in $\infty$, because $\bar\iota_{\infty}\cdot det(  h )= \iota_{\infty}$
and $\bar\iota_{\infty}$, $\iota_{\infty}$ are surjective
morphisms, and hence one deduces   the equalities (*).

Using   remark \ref{dt}, since
$$
 h_A  (\tau   )\wedge  \cdots\wedge  h_A  ( \tau^{n-1}     )\wedge
  [h_A \cdot \tau-\bar \tau\cdot h_A ](\tau^{n-1}  )   $$
is an element of $R\{\tau_{det}\}$ of degree $\leq nd  $, we have
that
$$deg_{\bar \tau}[h_A \cdot \tau-\bar \tau\cdot h_A ](\tau^{n-1} ) \leq n d+n-1,$$
and we conclude that $[h_A \cdot \tau-\bar \tau\cdot h_A
](\tau^{n-1} )=0$ because
$$[h_A \cdot \tau-\bar \tau\cdot h_A ](\tau^{n-1} )\in R\{\bar \tau\}\cdot I $$
and $deg(I)=d+1$. Thus, $h_A :R\{\tau\}\to R\{\bar\tau\}$ is an
isogeny of degree $nd$.

\end{proof}

\begin{lem} Let $X=Spec(A)$ be a smooth, noetherian scheme of  dimension $2n$. Let   $  Z_1 + \cdots + Z_r  $
be an $n$-cycle  in $X$   such that   $Z_i$ are different
irreducible closed subschemes  of dimension $n$ in $X$. If the
closed subscheme  $Z:= Z_1\cup  \cdots \cup Z_r $  is given by an
ideal generated by $n$ elements $a_1,\cdots,a_n\in A$, then  the
$n$-cycle  $  Z_1 + \cdots + Z_r   $   is rationally equivalent to
$0$.
\end{lem}
\begin{proof} Let ${\mathfrak I}$ be an ideal within $A$.  We denote by $Z({\mathfrak I}) $   the cycle associated in $X$
with the closed subscheme given by ${\mathfrak I}$.  The prime
ideal within $A$ given by $Z_i$ is denoted by $P_i$. Thus,
$$Z(P_1\cap \cdots \cap P_r)= Z_1 + \cdots + Z_r .$$

Let us consider  the ideal $(a_2,\cdots,a_n)$ within $A$ generated
by $ a_2,\cdots,a_n $ and let $ Q_1\cap \cdots \cap Q_h$ be a
minimal primary decomposition of this ideal.  If $Y_1 , \cdots ,
Y_k,(k\leq h)$ are the irreducible components of the closed
subscheme within $X$ given by $(a_2,\cdots,a_n)$, then dim$
Y_j\geq n+1$. We may assume, reordering the indices, that
$Z(Q_1)=Y_1,\cdots , Z(Q_k)=Y_k$.

By taking the localisation by $P_i$, one obtains
$$(A/Q_1\cap \cdots \cap Q_h )_{P_i}, $$
which  is a local ring  of dimension $1$ because $Y_1 , \cdots ,
Y_k$ has dim$ Y_j\geq n+1$. From the equality of rings
$$ A/Q_1\cap \cdots \cap Q_h+(a_1)=A/P_1\cap \cdots \cap P_r$$
  one obtains
$$(A/Q_1\cap \cdots \cap Q_h+(a_1))_{P_i}=(A/P_i)_{P_i}.$$
Therefore, $(A/Q_1\cap \cdots \cap Q_h )_{P_i} $      is principal
and hence an integral domain, and therefore there exists a unique
$Q_{l_i}, (l_i\leq k)$ with $Q_{l_i}\subset P_i$. If we denote by
$P_{j_1},\cdots, P_{j_{m_j} }$ the $P_i$'s with $Q_j\subset
P_{j_1},\cdots,Q_j\subset P_{j_{m_j} }, ( j\leq k) $, then within
the $n+1$-dimensional scheme, $Z(Q_j)=Y_j$
$$Z_{j_1}+\cdots  +Z_{j_{m_j} }$$ is given by the zero locus of $a_1$, which proves that $Z_{j_1}+\cdots
+Z_{j_{m_j} }$ is   rationally equivalent to $0$ on $X$.

We conclude because $  Z_1 + \cdots + Z_r
=\sum_{j=1}^kZ_{j_1}+\cdots  +Z_{j_{m_j} }$.

\end{proof}
\begin{thm}\label{ct} The correspondence
  $$T^n_I:=T(n\cdot d)+\Gamma(Fr)*T(n\cdot d-1)+\cdots +\Gamma(Fr^{n\cdot d-1})* T (1)+\Gamma(Fr^{n\cdot d })
$$
is trivial(= rationally equivalent to $0$ as an $n$-cycle within
${\mathcal E}_{n,\vert {\Bbb P}_1\vert_{nd}}^{I \infty}\times
{\mathcal E}_{n,\vert {\Bbb P}_1\vert_{nd}}^{I \infty}$).
\end{thm}
\begin{proof}
  Bearing in mind   corollary \ref{tri} and   lemma \ref{condiciones}, this correspondence is given by the zero locus of
$n$ regular functions of ${\mathcal E}_{n,\vert {\Bbb
P}_1\vert_{nd}}^{I \infty}\times {\mathcal E}_{n,\vert {\Bbb
P}_1\vert_{nd}}^{I \infty}$.

We are within the hypothesis of the above lemma because the
projection over the first entry, $T(r)\to {\mathcal E}_{n,\vert
{\Bbb P}_1\vert_{nd}}^{I \infty}$ is an \'{e}tale morphism and
therefore  $T(r)$ is smooth. Moreover, because of lemma \ref{un}
if $i\neq j $ then:
$$\Gamma(Fr^{i })* T(n\cdot d-i)\cap \Gamma(Fr^{j })* T(n\cdot d-j)=\emptyset,$$
with $  i,j\leq n\cdot d$.

\end{proof}

\subsection{ Some explicit calculations }
One can make explicit calculations by using the anti-equivalence
between elliptic sheaves and Drinfeld modules (c.f:\cite{ Dr2},
\cite{Mu}) and by   using the explicit calculation of the global
sections "$s$", (\cite{Al1}, remark 3.1), in terms of the
$I$-torsion elements of the Drinfeld modules. For $n=1$,
calculations are made in \cite{ An2} and in the spirit of this
work in \cite{Al3}: (example 2, page 21) and in \cite{Al2}, 3.2.

We begin with the following example.

\begin{ex} $n=1$, $I=p(t){\Bbb F}_q[t]$, $p(t)=(t-a_1)\cdots
(t-a_{d+1})$, with $ a_i\neq a_j$ for $i\neq j$ and $ a_i\in {\Bbb
F}_q$. Let $( L_j,i_j,\tau )$ be the rank-one elliptic sheaf
defined over the Carlitz's cyclotomic ring $K_1^{I \infty }={\Bbb
F}_q(\bar t)[\delta]$,
 with $\delta$ an element of a closed field of ${\Bbb F}_q(t)$
 verifying:
  $$  \phi_{p(t)}(\delta)= \delta^{q^{d+1}}+\cdots+c_1\delta^{q } +p(\bar t)\delta=0,$$
where $\phi$ is the Drinfeld module $\phi_t=\tau +\bar t$, (remark
3 of section 2.2).

Let us consider the $I\infty$-level structure, $ \iota_{I\infty}
$,  for $(  L_j,i_j,\tau )$. We have
$$\iota_{I }:L_0\to K_1^{I \infty }[t]/p(t)$$
given by $\iota_{I\infty}(s)=m_1  \delta_1 \frac{p(t)}{ t-\alpha_1
} + \cdots+m_{d+1} \delta_{d+1}  \frac{p(t)} {  t-\alpha_{d+1}}$
and $$\iota_{ \infty}:L_0\to K_1^{I \infty }[t^{-1}]/t^{-1},$$
with $\iota_{ \infty}(s)=1$. Here $L_0=s\cdot \o_{{\Bbb
P}^1}\otimes K_1^{I \infty }$, $ \phi_{\frac{p(t)}{
t-a_j}}(\delta)=\delta_{j}$ and the $m_j$ are obtained from the
equality
$$\frac{1}{p(t)}= \frac{m_1} { t-a_1} +\cdots+\frac{m_{d+1} } { t-a_{d+1}} .$$

We shall obtain the  element of $K_1^{I \infty }\otimes K_1^{I
\infty }$ whose divisor is
$$ T(   d)+\Gamma(Fr)*T(
d-1)+\cdots +\Gamma(Fr^{  d-1})* T (1)+\Gamma(Fr^d ).$$

Let $\pi_1$ and $\pi_2$ be the natural projections
$$Spec(K_1^{I \infty }\otimes K_1^{I \infty })\to Spec(K_1^{I \infty
}).$$ The morphism, $h $, of lemma \ref{morfismo} applied to the
rank-one line bundles  with a $I$-level structure,
$\pi_1^*(L_0,\iota_{I })$ and $\pi_2^*(L_0(d),\iota_{I })$, is
given by:
$$ h(\pi_1^*s)=[m_1 \frac{\delta_1\otimes 1}{1\otimes \delta_1} \frac{p(t)}{ t-\alpha_1 } +
\cdots+m_{d+1}\frac{\delta_{d+1}\otimes 1}{1\otimes
\delta_{d+1}}\frac{p(t)} {  t-\alpha_{d+1}}] \pi_2^*s.$$

By lemma \ref{condiciones}, one must impose on $h$ that:
$$\pi_2^*\iota_{
\infty}(h(\pi_1^*s))=\pi_1^*\iota_{I\infty}(\pi_1^*s)$$ by the
definition of $\infty$-level structures
$$\pi_2^*\iota_{
\infty}(h(\pi_1^*s))=m_1 \frac{\delta_1\otimes 1}{1\otimes
\delta_1}   + \cdots+m_{d+1}\frac{\delta_{d+1}\otimes 1}{1\otimes
\delta_{d+1} } $$ which is the leader coefficient of the
polynomial
$$m_1 \frac{\delta_1\otimes 1}{1\otimes \delta_1}  \frac{p(t)}{t-a_1} +
\cdots+m_{d+1}\frac{\delta_{d+1}\otimes 1}{1\otimes
\delta_{d+1}}\frac{p(t)}{ t-a_{d+1} }.$$

 Since
 $\iota_{I\infty}( s)=1$,
 the element sought   is
$$m_1 \frac{\delta_1\otimes 1}{1\otimes \delta_1}   +
\cdots+m_{d+1}\frac{\delta_{d+1}\otimes 1}{1\otimes \delta_{d+1}
}-1.$$ \end{ex}

\begin{ex}  Now, we follow with the easiest  anabelian case.

$n=2$, $I= t {\Bbb F}_q[t]$. Let $\phi_t=a\sigma^2+b\sigma+\bar t$
be a Drinfeld module of rank two defined over the ring
$$B_2^{I\infty }=({\Bbb F}_q[a,a^{1/1-q},b, \bar t, r(\bar
t)^{-1}]/a+b+\bar t-1)[\Gamma],$$ with $\phi_t(1)=0$,
$\Gamma^q-\Gamma\neq 0$ and $\phi_t(\Gamma)=0$.
 $r(t)$ is the product of the monic polynomials of degree less
 than
or equal to $2$. Let $(  M_j,i_j,\tau   )$ be the rank-two
elliptic sheaf associated with $\phi$, and  let $ \iota_{I\infty}
$ be an $I\infty$-level structure for $(  M_j,i_j,\tau   )$ given
by:
$$\iota_{I }:M_0\to (B_2^{I\infty }[t]/ t)^{\oplus 2}\simeq (B_2^{I\infty } )^{\oplus 2}, $$
given by $\iota_{I\infty}(s)=(1,\Gamma) $ and
$\iota_{I\infty}(\tau s)=(1,\Gamma^q)$. Recall that:
$$M_0=s\cdot
(\o_{{\Bbb P}^1}\otimes B_2^{I\infty }) \oplus \tau s\cdot
(\o_{{\Bbb P}^1}\otimes B_2^{I\infty }).$$

The $\infty$-level structure
$$\iota_{ \infty}:Det (M_0)\to B_2^{I\infty }[t^{-1}]/t^{-1}$$
is given by $\iota_{ \infty}(s\wedge \tau s)=a$.

Let $\pi_1$ and $\pi_2$ be the natural projections:
$$Spec(B_2^{I\infty }\otimes B_2^{I\infty })\to Spec(B_2^{I\infty }).$$
The morphism, $h$, of lemma \ref{morfismo} applied to the rank-two
vector bundles with a $I$-level structure, $\pi_1^*(M_0,\iota_{I
})$ and $\pi_2^*(M_0,\iota_{I })$, (here $d=0$) is given by the
matrix product
$$ D:=\left(\begin{matrix} 1   & 1   \\ 1\otimes\Gamma   & 1\otimes\Gamma^q
\end{matrix}\right)^{-1}\cdot  \left(\begin{matrix} 1   & 1   \\ \Gamma \otimes 1 & \Gamma^q \otimes 1\end{matrix}\right).$$

By lemma \ref{condiciones}, one must impose on $h$, that $deg
(h_A(1))=0$, where
$$h_A:B_2^{I\infty }\otimes B_2^{I\infty }\{\tau
\}\to B_2^{I\infty }\otimes B_2^{I\infty }\{\tau \}$$ is the
restriction of $h$ to $ {{\Bbb P}^1}\setminus \{\infty\}$. By
considering the second entry of $D\left(\begin{matrix} 1       \\
0
\end{matrix}\right)$, this condition is:
$$  *:= (1\otimes\Gamma^q   - 1\otimes\Gamma )^{-1}(\Gamma \otimes 1 -
1\otimes\Gamma)=0.$$

We must now impose that:
$$\pi_2^*\iota_{
\infty}(det(h)(\pi_1^*(s\wedge \tau
s)))=\pi_1^*\iota_{I\infty}(\pi_1^*(s\wedge \tau s)).$$
  Since,
 $\iota_{ \infty}(s\wedge \tau s)=a$, we have:
$$\pi_2^*\iota_{ \infty}(det(h)(\pi_1^*(s\wedge \tau s) ))=Det(D)\cdot
(1\otimes a),$$
 and
 $$\pi_1^*\iota_{I\infty}(\pi_1^*(s\wedge \tau
s))=a\otimes 1.$$  Thus we obtain the element:
$$**:= (\Gamma^q -\Gamma) \otimes  (\Gamma^q   -
 \Gamma)^{-1} - a\otimes a^{-1}   $$

The final result is that the diagonal   subscheme of
$Spec(B_2^{I\infty }\otimes B_2^{I\infty })$ is the zero locus of
the ideal generated by $*$ and $**$. \end{ex}

\section{ The additive case: $n=2$ (anhilators for cusp forms) }
In this section, we shall follow  the notation setout in the
introduction. The set of cusps is $\overline {{\mathcal
E}}({I\infty}) \setminus {{\mathcal E}({I\infty})}$.   We denote
by ${\mathcal C}^0({I\infty})$ the divisor class group on
$\overline {{\mathcal E}}({I\infty}) $ whose support lies among
the cusps. As in the introduction we follow  the notation and
results of \cite{GR}. For the definition and study of cusps forms,
readers are referred to the works of Gekeler, Goss or the
Habilatationshrift of Gebhard B\"ockle.

We now prove a Lemma which is the counterpart for $n=2$ for
Stickelberger's Theorem.
\begin{lem}  $T(2d   )+  T(2d -1)+\cdots +  T
(1)+\Gamma(Id) $ annihilates the group $Pic({\mathcal
E}({I\infty}))$.
\end{lem}
\begin{proof} This lemma is proved by using   theorem \ref{ct}, and the fact that the divisor group, ${\mathcal D}^0({\mathcal
E}({I\infty})$, of the affine curve, ${{\mathcal E}({I\infty})}$
defined over $K^\infty_I$  is a subgroup of the group of Weil
divisors of ${\mathcal E}_{2,\vert {\Bbb P}_1\vert_{2d}}^{I
\infty}$. Recall that  ${\mathcal E}_{2,\vert {\Bbb
P}_1\vert_{2d}}^{I \infty}$  is a smooth variety of dimension 2.
\end{proof}

 Note that the Hecke correspondences operate  over
the cusps. Thus, the above correspondence  gives a group
endomorphism, ${\mathcal C}^0({I\infty})\to  {\mathcal
C}^0({I\infty})$. We denote by  $\overline {S}_2(d)$, $  {S}_2(d)$
and $ {S'}_2(d)$ the group endomorphisms given by,
$$T(2d )+
T(2d  -1)+\cdots + T (1)+\Gamma(Id) $$ over the groups
$Pic^0(\overline {{\mathcal E}}({I\infty}) )$, $Pic({\mathcal
E}({I\infty}))$  and ${\mathcal C}^0({I\infty})$,
  respectively.

Let us consider  $j^*$,   the pull back of the line bundles over
$\overline {\mathcal E}({I\infty}) $ by the natural inclusion:
$$j:  {{\mathcal E}({I\infty})}\hookrightarrow \overline {\mathcal
E}({I\infty}) .$$ We assume that:
$$j^*:Pic^0(\overline {\mathcal E}({I\infty}) )\to Pic({\mathcal
E}({I\infty}))$$ is surjective. For example if in the   cusps
points  $\overline {{\mathcal E}}({I\infty}) \setminus {{\mathcal
E}({I\infty})}$  there exists a rational point over $K_I^\infty$.
If this does not occur then it suffices to change
  $Pic({\mathcal E}({I\infty}))$ by $Pic^0({\mathcal E}({I\infty}))$.

\begin{lem}  If $Pic({\mathcal
E}({I\infty}))$ is an infinity group, then $Ker ( \overline
{S}_2(d))$ also has infinity elements.
\end{lem}
\begin{proof} If $Coker(S)$ is
not finite then we finish the proof, since   ${\mathcal
C}^0({I\infty})$ is a finite generated group. Thus, we can assume
that $ ker({S'}_2(d))$ and $coker({S'}_2(d))$ are finite groups.

 From the serpent lemma applied to the commutative diagram
$$ \xymatrix {  0\ar[r]&{\mathcal C}^0({I\infty})\ar[d]^{S'_2(d)}\ar[r]&
Pic^0(\overline {\mathcal E}({I\infty})  )
\ar[d]^{\overline{S}_2(d)}\ar[r]^{j^*} &  Pic({\mathcal
E}({I\infty})) \ar[d]^{ {S}_2(d)}\ar[r]^{j} \ar[r]&0
\\0 \ar[r]&{\mathcal C}^0({I\infty}) \ar[r]&Pic^0(\overline {\mathcal E}({I\infty}))  \ar[r]^{j^*} &  Pic({\mathcal
E}({I\infty}))\ar[r] &0}
$$
one obtains an exact sequence:
$$Ker (S'_2(d))\to Ker (\overline{S}_2(d)) \to Pic({\mathcal
E}({I\infty})) \overset {\delta} \to Coker (S'_2(d)).$$ We
conclude since   $Ker (S'_2(d))$, and  $Coker (S'_2(d))$ are
finite groups and because by hypothesis,  $Pic({\mathcal
E}({I\infty}))$  is an infinity group.
\end{proof}

\begin{thm}  If the the cardinal of the group $Pic ({\mathcal
E}({I\infty}))$ is $\infty$, then there exists a cusp form for
$\Gamma_{I\infty}$ that is annihilated by $\tilde T (2d )+ \tilde
T (2d -1)+\cdots + \tilde T (1)+ Id $.
\end{thm}
\begin{proof} We denote by  $J$ the Jacobian of the curve $\overline {\mathcal
E}({I\infty})$ over $K_I^\infty$. Thus, the correspondence
$\overline{S}_2(d)$ gives an endomorphism of this Jacobian. By the
last lemma, this endomorphism can not be an isogeny. Accordingly,
  the morphism induced  over the tangent space  $T_e(J)$   of $J$
over the zero element,
$$\tilde{S}_2(d):T_e(J)\to T_e(J)$$ is not
injective. We conclude because  the tangent space, $T_e(J)$, is
the dual  space of the space of $1$-holomorphic differential
forms, $ H^0(\overline {\mathcal E}({I\infty}), \Omega_{\overline
{\mathcal E}({I\infty})}) $, and the space of cusp forms is
identified with
$$H^0(\overline {M}_{ \Gamma_{I\infty} }, \Omega_{\overline {M}_{ \Gamma_{I\infty} }})=
H^0(\overline {\mathcal E}({I\infty}), \Omega_{\overline {\mathcal
E}({I\infty})})\otimes _{K_I^\infty}C.$$
\end{proof}

\section{Ideal class group anhilators for the cyclotomic function fields}

 We consider  ${\mathcal E}_{1,\vert {\Bbb P}_1\vert_{nd}}^{I \infty}=Spec( {\mathcal B}_1^{I
\infty})$. The construction of ${\mathcal B}_1^{I\infty}$ is
detailed in section 2.2, remark 3, and is essentially as follows:
Let $((L_j,i_j,\tau),\iota_{ I\infty})$ be an element of
${\mathcal E}_{1,\vert {\Bbb P}_1\vert_{nd}}^{I \infty}$.  To
construct ${\mathcal B}_1^{I\infty}$,  we can fix a global section
$s$ of $L_0$ such that $t\cdot s=\lambda\cdot s + \tau s$, and
$\iota_\infty (s)=1$. Hence,
  $ Spec({\mathcal B}_1^{I \infty})$ represents the pairs  $(\phi, \iota_I)$,  with $\phi$ a rank $1$-normalized
Drinfeld module  and $\iota_I$ an $I$-level structure for $\phi$.
${\mathcal B}_1^{I \infty}$    is considered in \cite{An1} and
\cite{C} and  is obtained from the $I$-torsion elements of a
normalized   Drinfeld module.  The "zero" morphism ${\mathcal
E}_{1,\vert {\Bbb P}_1\vert_{nd}}^{I \infty}\to Spec({\Bbb
F}_q[t])$ gives a   Galois extension $K^\infty_I/{\Bbb F}_{q}(t)$
of group $({\Bbb F}_q[t]/I)^\times$. We denote by $Y^\infty_I$ the
proper smooth curve associated with ${\mathcal E}_{1,\vert {\Bbb
P}_1\vert_{nd}}^{I \infty}$.

We consider the Hecke correspondence:
$$T(J_1 ,\cdots , J_n)\subset {\mathcal E}_{n,\vert {\Bbb P}_1\vert_{nd}}^{I
\infty}\times {\mathcal E}_{n,\vert {\Bbb P}_1\vert_{nd}}^{I
\infty},$$ which is of degree $d(J_1,\cdots, J_n)$    over the
second component. Let $J$ be the product of ideals:
$$\prod^n_{i=1} J_i:=J.   $$
  $T(J)$ denotes the Hecke correspondence on ${\mathcal E}_{1,\vert {\Bbb P}_1\vert_{nd}}^{I
\infty}$ given by $J$.

There exists an action, $T(J_1,\cdots, J_n)^*$ and $T(J)^*$, of
these    correspondences over the functors, $\underline {Pic}(
{\mathcal E}_{n,\vert {\Bbb P}_1\vert_{nd}}^{I \infty} )$ and $
\underline {Pic}^0( Y^\infty_I)$, respectively. These functors are
defined over the category of $ {\Bbb F}_q$-schemes. Recall that
the projections $\bar \pi_1$ and $\bar \pi_2$
$$T(J_1,\cdots, J_n)\to   {\mathcal E}_{n,\vert {\Bbb
P}_1\vert_{nd}}^{I\infty}$$

are \'{e}tale. In this way it is possible to define $T(J_1,\cdots,
J_n)^*(:=(\bar\pi_2)_*\cdot\bar\pi_1^*$).

\begin{rem}
Let us consider   the morphism $det:{\mathcal E}_{n,\vert {\Bbb
P}_1\vert_{nd}}^{I \infty} \to  {\mathcal E}_{1,\vert {\Bbb
P}_1\vert_{nd}}^{I \infty}$. We set
$$T(J_1,\cdots, J_n)^* {det}^*[  D]=  d(J_1,\cdots, J_n)   {det}^* T(J)^*[  D] $$

    $[  D]$ is the class of a divisor $  D$ on $  {\mathcal E}_{1,\vert {\Bbb
P}_1\vert_{nd}}^{I \infty}  $.

This result is proved bearing in mind the projection formula for
the rational equivalence of cycles; that  $\bar \pi_2$ is an \'{e}tale
morphism of degree $d(J_1,\cdots, J_n)$, and that  given $e:=
(E,\iota_{I\infty}) \in  {\mathcal E}_n^{I\infty}$  we have:
$$\bar \pi_2[T(J_1,\cdots, J_n)\cap (  det^{-1}(det(e ) ) \times {\mathcal E}_{n,\vert {\Bbb P}_1\vert_{nd}}^{I\infty})]=
 det^{-1}(  T(J)^* det(e)),$$
 $$T(J_1,\cdots, J_n)\cap (   det^{-1}(det(e ) ) \times {\mathcal E}_{n,\vert {\Bbb P}_1\vert_{nd}}^{I\infty})
=\bar\pi_2^{-1}[
 det^{-1}( T(J)^* det(e))].$$

\end{rem}

\begin{lem} The correspondence
$$ D_I^n :=\sum_{i=0}^{nd-i} \Gamma (Frob^i) * [\sum_{ \underset{ J+ I={\Bbb F}_q[t] , deg(J)=i}
{J\subset {\Bbb F}_q[t]  }}  (\sum_{\underset{\prod^n_{k=1} J_k
=J} {J_1\subseteq \cdots \subseteq J_n}}  d(J_1,\cdots, J_n))
T(J)]$$ is trivial on $Y^\infty_I\times Y^\infty_I$ up to vertical
and horizontal correspondences.
\end{lem}
 \begin{proof} It suffices to consider a curve $ Z\to {\mathcal E}_{n,\vert {\Bbb
P}_1\vert_{nd}}^{I\infty}$ such that the morphism composition:
$$g:Z\to {\mathcal E}_{n,\vert {\Bbb
P}_1\vert_{nd}}^{I\infty}\overset{det}\to {\mathcal E}_{1,\vert
{\Bbb P}_1\vert_{nd}}^{I\infty}$$ is not constant.  By the above
remark, $(det )^*\cdot (D_I^n)^*=(T_I^n)^*\cdot (det )^*$. Since
$T_I^n$ is rationally equivalent to zero, we have that $(g )_*(g
)^*\cdot (D_I^n)^*$ is trivial on $Y^\infty_I\times Y^\infty_I$ up
to vertical and horizontal correspondences, but  by the projection
formula
$$(g )_*(g )^*\cdot (D_I^n)^*=m\cdot D_I^n,$$
with $m\in{\Bbb N}$. We conclude, since   the ring of
correspondences, modulo horizontal and vertical ones, is without
${\Bbb Z}$-torsion.

\end{proof}
By using the Euler  product of lemma \ref{euler} one could find
another proof of this lemma.
    \bigskip

 We consider $J =q(t)  {\Bbb F}_{q}[t]  $  with $q(t)$ monic. $T(J)$  is given by the graphic, $\Gamma (q(t))$, of
the automorphism of ${\mathcal E}_{1,\vert {\Bbb
P}_1\vert_{nd}}^{I \infty}$ obtained from the action of $ q (t)\in
({\Bbb F}_q[t]/I)^\times  $. Recall that to obtain ${\mathcal
B}_1^{I \infty}$ we have fixed a global section $s$ of $L_0$ such
that $t\cdot s=\lambda\cdot s + \tau s$, and $\iota_\infty (s)=1$.
In this way,  $T(J)=\Gamma (q(t) ) $.  By section 2.3, if we set
$J_i=q_i(t)  {\Bbb F}_{q}[t]  $, then:
$$\varphi(q(t),n):=\sum_{\underset{\prod^n_{k=1} J_k =J} {J_1\subseteq \cdots \subseteq J_n}}  d(J_1,\cdots, J_n))$$
is the   number of submodules $N\subseteq {\Bbb F}_q[t]^{\oplus
n}$ such that:
$$ {\Bbb F}_q[t]^{\oplus n}/N\simeq {\Bbb F}_q[t]/q_1(t)\oplus \cdots \oplus {\Bbb F}_q[t]/q_n(t),$$
with the product of the invariant factors $q_1(t)  \cdots  q_r(t)$
equal to $q(t)$.  Therefore, if we consider $p(t){\Bbb F}_q[t]=I$,
then:

\begin{cor}The correspondence
$$  \sum_{i=0}^{nd} [\Gamma (Frob^{nd-i})  *(\sum_{ \underset{(p(t),q(t))=1 , deg(q(t))=i} {q(t) \text{ monic}\in
{\Bbb F}_q[t]  }}    \varphi(q(t),n)\cdot \Gamma(q(t)  )]$$ is
trivial on $Y^\infty_I\times Y^\infty_I$ up to vertical and
horizontal correspondences.
\end{cor}
 \bigskip
 \begin{ex} We can check this result for $n=2$ and $p(t)=t(t-1)$. Let $K^\infty_{t(t-1)}/{\Bbb
F}_{q}(t)$ be    the Galois extension  of group $({\Bbb
F}_q[t]/t(t-1))^\times$.

One has that $\varphi(t-\alpha,2)=q+1$,
$\varphi((t-\alpha)^2,2)=q^2+q+1$,
 $\varphi((t-\alpha)(t-\beta) ,2)=q^2+2q+1$, and $\varphi(t^2+at+b,2)=q^2+1$ with $t^2+at+b\in{\Bbb F}_q[t]$ an irreducible
polynomial and $\alpha, \beta \in {\Bbb F}_q$, $\alpha \neq
\beta$. Thus
$$ (*)\quad \sum_{i=0}^{2} [\Gamma (Frob^{2-i})  *(\sum_{ \underset{(t(t-1),q(t))=1 , deg(q(t))=i} {q(t) \text{ monic}\in
{\Bbb F}_q[t]  }}    \varphi(q(t),2)\cdot \Gamma(q(t)  )] $$ is
$$   \Gamma (Frob^2 )  +  (q+1) \sum_{   \alpha\neq 0,1  }
\Gamma (Frob  )*\Gamma(t-\alpha  ) +$$
$$+(q^2+2q+1) \sum_{\underset{  \underset {\alpha\neq\beta} {\alpha,\beta \neq 0,1} }  { \{\alpha ,\beta \}\subset {\Bbb
F}_q } } \Gamma( (t-\alpha)(t-\beta)) +$$
$$+(q^2+ q+1 ) \sum_{  \alpha\neq 0,1  }
\Gamma( (t-\alpha)^2)+$$
$$+(q^2+1 ) \sum_{   \underset {t^2+at+b, \text{irreducible}} {a,b\in
{\Bbb F}_q  }   } \Gamma  (t^2+at+b),$$ and grouping terms, we
have:
$$  \Gamma (Frob )* [\Gamma (Frob )  +  \sum_{    \alpha\neq 0,1  }
\Gamma(t-\alpha  )]  +$$
$$+q[ \sum_{    \alpha\neq 0,1   }
\Gamma(t-\alpha  )] *[\Gamma (Frob )  +  \sum_{    \alpha\neq 0,1
} \Gamma(t-\alpha  )]+$$
$$+(q^2+1)[\sum_{\underset{  \underset {\alpha\neq\beta} {\alpha,\beta \neq 0,1} }  { \{\alpha ,\beta \}\subset {\Bbb
F}_q } } \Gamma( (t-\alpha)(t-\beta))+ \sum_{\alpha\neq 0,1}
\Gamma( (t-\alpha)^2)+\sum \Gamma  (t^2+at+b)].$$ Now, bearing in
mind that the last sumand is:
$$(q^2+1)(\sum_{g\in ({\Bbb F}_q[t]/t(t-1))^\times} \Gamma(g)),$$
which is a trivial correspondence, we conclude that (*) is also
trivial because   the correspondence
$$   \Gamma (Frob )  +  \sum_{    \alpha\in
{\Bbb F}_q \setminus \{0,1\}   } \Gamma(t-\alpha  ) $$ is trivial
on $K^\infty_{t(t-1)}\otimes K^\infty_{t(t-1)}$ c.f. \cite{C}.
\end{ex}

\bigskip
\section{The above results without $\infty$ level structures}
With minor changes in the above results one can obtain   similar
results but over the modular varieties, ${\mathcal E}_{n}^{I  }$.
The results obtained match, for $n=1$,  the classical
Stickelberger's theorem over ${\Bbb Z}$. (c.f  \cite{Gr1},
\cite{Gr2}).

To obtain these results it suffices to replace in lemma
\ref{condiciones}  the condition imposed on $h $ to be a morphism
of $\infty$-level structures,  by the condition:
$$deg_{\bar
\tau}(h_A   (\tau^{n })-\bar \tau\cdot  h_A(
 \tau^{n-1})\leq
n-1+nd.$$ And now in corollary \ref{tri} one allows pairs,
$[(E,\iota_{I  } ), (\bar E,\bar\iota_{I  } )]$,  given by an
isogeny for $I$-level structures, $\Phi:(E,\iota_{I } )\to (\bar
E,\bar\iota_{I } )$, such that $\infty $ can be within $supp(\bar
E/\Phi (E))$. Thus, one obtains:
\begin{thm}\label{ct} The correspondence

 $  T(n\cdot d)+[\Gamma(Fr)   +\Gamma(Id)]*T(n\cdot d-1)+\cdots
    + [\Gamma(Fr^{n\cdot d-1}) +\cdots +\Gamma(Id)]*T (1) +$
    \newline$+
      [ \Gamma(Fr^{n\cdot d})  +\Gamma(Fr^{n\cdot d-1}) +\cdots+\Gamma(Fr)+\Gamma(Id)]
     $

  is trivial(= rationally equivalent to $0$ as an $n$-cycle) within
${\mathcal E}_{n,\vert {\Bbb P}_1\vert_{nd}}^{I  }\times {\mathcal
E}_{n,\vert {\Bbb P}_1\vert_{nd}}^{I  }  $.
\end{thm}

From the last theorem one has, for $n=2$:
\begin{lem}  The correspondence $T(2  d)+  2T(2 d-1)+\cdots +  2  dT
(1)+(2d+1)\Gamma(Id)) $ annihilates the group $Pic ({\mathcal
E}({I }))$.
\end{lem}
and

\begin{thm}  If the the cardinal of the group $Pic ({\mathcal
E}({I }))$ is $\infty$, then there exists a cusp form for
 $\Gamma_{I }$ that is annihilated by $$\tilde T(2  d)+ 2\tilde T(2
d-1)+\cdots +  2  d\tilde T (1)+(2d+1)\Gamma(Id)) .$$
\end{thm}

 {\bf{Acknowledgments} }I would like to
thank the referee for suggesting to deal with  annihilators for
cusp forms. I would like to thank N. Skinner for doing his best to
supervise my deficient English. I am also deeply grateful to
Ricardo Alonso Blanco and Jesus Mu\~noz Diaz for their help.



\vskip1.5truecm { \'Alvarez V\'azquez, Arturo}\newline {\it
e-mail: } aalvarez@gugu.usal.es

\end{document}